\documentclass[11pt,a4paper,twoside]{article}
\usepackage{amsfonts} 
\usepackage{amssymb}
\usepackage{amsmath, amsthm, array, ifthen}
\usepackage{mathrsfs} 
\usepackage{caption} 
\usepackage{latexsym}
\usepackage{enumerate}

\usepackage{enumitem}

\usepackage{xfrac}
\usepackage{graphicx}

\usepackage{lmodern}

\usepackage[french]{babel} 
\usepackage[utf8]{inputenc} 

\usepackage{color}

\usepackage[colorlinks]{hyperref}

\newcommand{\CC}{\mathbb{C}}

\newcommand{\epsi}{\varepsilon}

\newcommand{\mell}[2][f]{\mathcal{M}(f)} 

\newcommand{\diff}{\mathop{}\mathopen{}\mathrm{d}}

\newtheorem{prop}{Proposition}
\newtheorem{prop/not}{Proposition/notation}

\newtheorem{cor}{Corollaire}
\newtheorem{lem}{Lemme}

\theoremstyle{definition}
\newtheorem{rem}{Remarque}

\newtheorem{thmintro}{Théorème}

\newcommand{\ioe}{\leqslant}
\newcommand{\soe}{\geqslant}

\usepackage[usenames,dvipsnames,svgnames,table]{xcolor}

\title{ Conversions explicites des nombres premiers vers la fonction de M\"obius}
\usepackage[left=3cm,right=3cm,top=1.5cm,bottom=1.5cm]{geometry}
\author{Daval Florian}
\date{}
\begin{document}
\maketitle
\begin{abstract} Nous démontrons que $|\sum_{n \ioe x} \mu(n) | \ioe
x/160\,383$ pour tout $x\soe 8.4 \times 10^9$ améliorant ainsi le record pour ce type d'estimations qui était  $|\sum_{n \ioe x} \mu(n) | \ioe
x/4345$ pour tout $x\soe 2\,160\,535$ et était dû à Cohen, Dress et El Marraki. \\
Nous prouvons également que $|\sum_{n \ioe x} \mu(n) | \ioe
x/180\,194$ pour tout $x\soe  10^{19}$.
\end{abstract}

\section{Introduction}
Nous  utilisons  des résultats effectifs obtenus sur la fonction de Tchebichef $\psi(x)=\sum_{p^k \ioe x} \log p$ (où la somme porte sur toutes les puissances des nombres premiers inférieurs ou égaux à $x$), pour en obtenir pour la fonction sommatoire de la fonction de M\"obius (cette fonction arithmétique  notée $\mu(n)$ vaut $1$ ou $-1$ selon que $n$ est le produit d'un nombre pair ou impair de nombres
premiers différents, et $\mu(n)=0$  quand $n$ a des facteurs carrés). Les résultats actuels sur $\sum_{p^k \ioe x} \log p$  sont souvent meilleurs que pour $\sum_{n \ioe x} \mu(n)$,  il est donc avantageux d'essayer de convertir des estimations des nombres premiers vers la fonction de M\"obius. Nous obtenons dans ce papier le résultat suivant.
\begin{thmintro} \label{thm A}
\begin{equation*}
|\sum_{n \ioe x} \mu(n) | \ioe
\frac{x}{160\,383} \quad (x \soe 8.4 \times 10^9)\,.
\end{equation*}
\end{thmintro}
Donnons quelques points historiques. On sait depuis Landau que le théorème des nombres premiers peut se présenter sous plusieurs formes, en particulier les trois résultats suivants sont élémentairement équivalents : \begin{equation*}
M(x)=\sum_{n \ioe x} \mu(n)=o(x), \quad 
xm(x)=x \sum_{n \ioe x} \frac{\mu(n)}{n}=o(x), \quad 
\psi(x)-x=o(x) \;. 
\end{equation*}
On peut d'abord essayer de démontrer des versions effectives de ces résultats en remplaçant les $o(x)$ par des $O(x)$. Le premier pas dans ce sens a été accompli pour la fonction $\psi$ en 1850 par Tchebichef en utilisant la fonction $t \mapsto \lfloor t \rfloor- \lfloor t/2 \rfloor- \lfloor t/3 \rfloor- \lfloor t/5 \rfloor+ \lfloor t/30 \rfloor$.  Plus tard von Sternek adaptera la méthode pour la fonction $M$ avec la fonction
  $t \mapsto \lfloor t \rfloor-\lfloor t/2 \rfloor-\lfloor t/3 \rfloor-\lfloor t/6\rfloor$. Il  utilise en particulier l'inversion de M\"obius et von Sterneck arrive à prouver ainsi dans une
  publication de 1898 que $|M(x)| \ioe x/9+8$ pour tout $x\soe 1$ et que $|L(x)| \ioe x/9+\sqrt{x/2}+7$ où $L$ est la fonction sommatoire de la fonction de Liouville. 
  Ces méthodes élémentaires seront poursuivies par Costa Pereira dans~\cite{CostaPPsiM} pour $\psi$ et $M$ en utilisant une autre combinaison linéaire de dilatées de la fonction partie entière. Pour la fonction $\psi$ des résultats meilleurs seront obtenus par l'étude des zéros de la fonction zêta et une partie de ses résultats sera transférée par Schoenfeld vers la fonction $M$ dans l'article \cite{SchoenfeldMobius}, c'est  le premier à obtenir des majorations explicites de $|M(x)|/x$ qui tendent vers $0$. Pour $M(x)=O(x)$ le meilleur résultat est obtenu par des moyens élémentaires (c'est-à-dire sans utiliser de  zéros de la fonction zêta et sans passer par la fonction $\psi$)  dans l'article~\cite{MobiusCohenDressMarraki} (prépublié 11 ans avant), les auteurs démontrent que
 \begin{equation}
\text{(Cohen-Dress-El Marraki, 1996)} \quad  : \quad   |M(x)|  \ioe \frac{x}{4345}  \quad (x \soe  2\,160\,535)\,. \label{CDM}
\end{equation}
Nous utiliserons une modification de la méthode de Schoenfeld – déjà reprise par El
Marraki dans \cite{MarrakiMobiusIII} – en nous basant notamment sur une estimation effective de obtenue
par voie analytique par Faber et Kadiri dans \cite{faber2015new}. Elles démontrent que
\begin{equation*}
    |\psi(x)-x|\ioe 8.6386 \times 10^{-8}x  \quad (x \soe \exp(40)\approx 2.35 \times 10^{17})\,.
\end{equation*}
Pour obtenir notre amélioration nous modifions l'identité utilisée par  Schoenfeld mais il y a une autre différence, nous utilisons l'identité pour passer des versions effectives de $M(x)=O(x)$ et $\psi(x)-x=O(x)$ à une meilleure version effective pour $M(x)=O(x)$. Alors que Schoenfeld et El Marraki passaient de versions effectives de $M(x)=O(x \log^{-\alpha} x)$ et $\psi(x)-x=O(x \log^{-\beta} x)$ avec la condition $\alpha< \beta$, à une meilleure version effective pour $M(x)=O(x \log^{-\gamma} x)$ avec $\gamma=(\alpha+1)\beta/(\beta+1)$ proche de $\beta$ mais inférieur  (en itérant si possible). Ce qui ne laisse pas la possiblité de prendre $\alpha=\beta=0$. Par exemple  pour $\alpha=0$ et $\beta=2$, Schoenfeld part de
\begin{enumerate}
\item $|M(u)|/u \ioe 1/80$ pour tout $u\soe 1119$,
\item $|\psi(x)-\lfloor x \rfloor| \ioe 652 x/\log^2x$ pour tout $x>1$,
\end{enumerate}
et obtient 
\begin{equation*} 
\frac{|M(x)|}{x} \ioe \frac{1.2}{(\log x)^{2/3}}  \quad (x \soe \exp(900))\,.
\end{equation*}
Ensuite il reprend avec ce $\alpha=2/3$ qu'il vient d'obtenir et $\beta=2$ ce qui donne $\delta=10/9$,
\begin{equation*} 
\frac{|M(x)|}{x} \ioe \frac{26}{(\log x)^{10/9}}  \quad (x \soe \exp(1000))\,.
\end{equation*}
Si l'on remplace l'estimation pour $M$ numérotée 1 ci-dessus (qui est due à Macleod) par l'estimation \eqref{CDM} on obtient des meilleurs résultats. D'après Dress (cf. \S 2 de \cite{dress1977majorations} dernière phrase) dans les résultats effectifs obtenus par la méthode de Schoenfeld de la forme
\[
|M(x)| \ioe \frac{c_{\alpha} x}{(\log x)^{\alpha}} \quad (x \soe T)
\]
la constante $c_{\alpha}$ dépend en puissance $2/3$ environ  du $\epsi$ des majorations effectives du type $|M(x)| \ioe \epsi x$. Ainsi avec notre nouveau résultat (théorème \ref{thm A}) on peut s'attendre à une amélioration d'un facteur 10 par rapport à l'utilisation du résultat \eqref{CDM}. De plus l'estimation du théorème \ref{thm A} permet de ramener une version effective de $M(x)=O(x \log^{-1} x)$ qui sera valide pour des $x$ très grands par rapport à $10^{16}$, valeur en dessous de laquelle on dispose d'un modèle calculé 
directement pour $M$ (voir le lemme \ref{racine} page~\pageref{racine}) et qu'on se propose d'établir dans un prochain travail.
Pour obtenir ses résultats, Schoenfeld  passe d'abord par la fonction $N(x)=\sum_{n \ioe x}\mu(n) \log n$  et utilise le principe de l'hyperbole de Dirichlet  (equation (19) p. 227 dans \cite{SchoenfeldMobius}),
\begin{equation}
    -N(x)-1=\sum_{k \ioe y} (\Lambda(k)-1) M(x/k)+\sum_{j \ioe x/y} \mu(j) (\psi(x/j)-\lfloor x/j \rfloor)-(\psi(y)-\lfloor y \rfloor ) M(x/y) \label{identite Schoen}
\end{equation}
 qui est à comparer (à une division par $x$ près) à notre modification présentée dans le lemme \ref{identite} ci-après.
 
 Le travail \cite{from} de Ramaré fournit la meilleure estimation explicite qui tend vers $0$ pour $|M|$ à ce jour :
 \begin{equation} \label{Ramare M}
\frac{|M(x)|}{x} \ioe  \frac{0.130}{\log x}-\frac{0.118}{\log^2 x}  \quad (x \soe 1\,078\,853).
\end{equation}
Remarquons que l'estimation du théorème \ref{thm A} est meilleure jusqu'à environ $10^{9000}$, c'est comme on l'a dit une de ses utilités. On peut lire les grandes lignes de la preuve de Ramaré suivant un schéma analogue à celui de Schoenfeld avec une identité plus efficace. On pose $N_2(x)=\sum_{n \ioe x} \mu(n) \log^2 n$, $g(\ell)=\Lambda*\Lambda(\ell)-\Lambda(l)\log \ell +2\gamma$ et $G(x)=\sum_{\ell \ioe x}g(\ell)$, alors on a l'identité \begin{align*}N_2(x)+2\gamma&=\sum_{nj \ioe x} \mu(j) g(n)\\
&=\sum_{n \ioe y} g(n) M(x/n)+\sum_{n \ioe x/y} \mu(n) G(x/n)-G(y)M(x/y),\end{align*}
mais nous n'avons pas réussi pour l'instant à adapter ce travail pour $M(x)=O(x)$ avec le même succès. 
Nous allons pour finir présenter d'autres perspectives que l'amélioration de \eqref{Ramare M}. D'après nos travaux \cite{daval2019identites} et \cite{daval2020conversions} le résultat obtenu pour $M$ dans le théorème \ref{thm A} peut se convertir vers les fonctions $m(x)=\sum_{n\ioe x} \mu(n)/n$, $m_1(x)=m(x)-M(x)/x$ et  $\check{m}(x)=\sum_{n \ioe x}(\mu(n)/n) \log(x/n)$. La version lissée
$\check{m}(x)$ constitue un enjeu important de la 
théorie explicite des nombres premiers, comme en témoigne leur usage dans le travail  
d'Helfgott \cite{helfgott2012minor}  sur le problème de Goldbach ternaire. On pourra également consulter la récente publication de Zuniga Alterman \cite{alterman2022logarithmic} qui utilise des estimations explicites de cette fonction. On obtiendra, sauf problèmes nouveaux,
\begin{equation*}
|m_1(x)| \quad \text{ et } \quad |\check{m}(x)-1| \ioe   \frac{\pi^2}{6} \frac{1}{4345}\frac{1}{160\,383}+\epsi \overset{?}{ \ioe } 2 \times 10^{-9}  \quad ( x \soe  10^{16})\, 
\end{equation*}
en utilisant une fonction qui est au c\oe ur de l'article de Cohen, Dress et El Marraki \cite{MobiusCohenDressMarraki}, ce qui explique la présence de $1/4345$ et l'importance qu'il y aurait à poursuivre leurs idées. Dans cet article les trois auteurs donnent également des estimations de $Q(x)=\sum_{n \ioe x} \mu^2(n)$ qui sont interdépendantes de l'estimation sur $M$. D'ailleurs ces estimations sur $Q$  nous ont été utiles pour obtenir notre résultat.

\subsection*{Notations}
Autour de la fonction de Möbius $\mu(n)$, pour $x>0$ on note :
  \begin{itemize}
  \item $M(x)=\sum_{n \ioe x} \mu(n)$
  \item $m_1(x)=\sum_{n\ioe x} \tfrac{\mu(n)}{n}(1-n/x)=m(x)-M(x)/x$
  \item $N(x)=\sum_{n \ioe x}\mu(n) \log n$\;.
  \end{itemize}
  
Autour de la fonction de von Mangoldt $\Lambda(n)$, pour $x>0$ on note :
  \begin{itemize}
  \item $\psi(x)=\sum_{n \ioe x} \Lambda(n)$\;.
  \end{itemize}
  
  Soit $R\soe 0$ une fonction définie sur un intervalle $I$.

  \begin{itemize}
  \item $f(x)=O(R(x))$ pour $x\in I$ signifie qu'il existe un nombre $c$ tel que $|f(x)| \ioe c R(x)$
    pour tout $x \in I$.

  \item $f(x)=O^*(R(x))$ pour $x\in I$ signifie que $|f(x)| \ioe  R(x)$
    pour tout $x \in I$.
  \end{itemize}
\section{Principe général}
Nous présentons et prouvons dans cette section le lemme \ref{identite} ci-dessous qui est le moteur principal de notre démonstration. 

Une des clefs de l'amélioration de cet article est de ne pas passer par $\Lambda-1$ bien que la convolution entre $1$ et $\mu$ soit intéressante. Cela  fait gagner un facteur $2$ dans la majoration de la première somme du lemme \ref{identite}. En effet, dans le travail de Schoenfeld \cite{SchoenfeldMobius}  (voir l'identité \eqref{identite Schoen} ci-dessus) et repris par El Marraki dans \cite{MarrakiMobiusIII}, apparaît la fonction arithmétique $ |\Lambda(k)-1|$  qui a une densité naturelle de $2$, dans notre adaptation on a  $ |\Lambda(k)|$ qui a pour densité $1$.  En contrepartie la fonction $m_1$ apparaît mais dans notre thèse \cite{daval2019identites} nous avons obtenu de bonnes estimations de cette fonction.

Puisque nous visons des estimations sur la fonction $M(x)$ il faut avoir en tête que la fonction $N(x)=\sum_{n \ioe x}\mu(n) \log n$  est très proche de $M(x)\log x= \sum_{n \ioe x}\mu(n) \log x$.
\begin{lem} \label{identite}
Pour tous nombres réels $x$ et $y$ avec $1 \ioe y \ioe x$ on a l'égalité :
\begin{equation*}
    \frac{-N(x)}{x}=\sum_{k \ioe y} \frac{\Lambda(k)}{k} \frac{M(x/k)}{x/k}\\
    +\sum_{j \ioe x/y} \frac{\mu(j)}{j} \frac{\psi(x/j)-x/j}{x/j}-\frac{\psi(y)-y}{y}\frac{M(x/y)}{x/y}+m_1(x/y)\;. 
\end{equation*}
\end{lem}

Nous n'utiliserons pas les mêmes démarches pour les valeurs de $x$ plus grandes que $10^{19}$ et celles avant que nous appellerons ici les petites valeurs. Dans tous les cas c'est le lemme \ref{identite} qui estimera la fonction $|N|$.

Pour les petites valeurs, nous utiliserons les modèles suivants pour $|M|$ et $|\psi|$. Pour $|M|$ le modèle est du à Hurst, il a calculé  les valeurs de $M(n)$ dans l'intervalle donné alors que le modèle pour $\psi$ a été obtenu analytiquement par B\"{u}the en utilisant en particulier des versions effectives de la formule explicite de von Mangoldt pour $\psi$ et des interpolations. 
\begin{lem}[Modèles en racine carrée] \label{racine} On a
\begin{align*}
|\psi(X)-X| \ioe 0.94 \sqrt{X} \quad &( 11 \ioe X \ioe 10^{19})\,, \\
|M(X)| \ioe 0.571 \sqrt{X} \quad &( 33 \ioe X \ioe 10^{16})\,.
\end{align*}
\end{lem}

\begin{proof}[Démonstration du lemme \ref{identite}]
L'identité entre séries génératrices $(1/\zeta)'=(1/\zeta).(-\zeta'/\zeta)$ donne la convolution entre fonctions arithmétiques $-\mu.\log=\mu \star \Lambda$ et donc 
\begin{equation*}
    -N(x)=\sum_{jk \ioe x} \mu(j) \Lambda(k) \quad (x\soe 1). 
\end{equation*}
Le principe de l'hyperbole de Dirichlet conduit à 
\begin{equation}
    -N(x)=\sum_{k \ioe y} \Lambda(k) M(x/k)+\sum_{j \ioe x/y} \mu(j) \psi(x/j)-\psi(y)M(x/y) 
\end{equation}
puis en faisant apparaître $\psi(t)-t$ plutôt que $\psi(t)$ on obtient que $-N(x)$ vaut
\begin{equation*}
    \sum_{k \ioe y} \Lambda(k) M(x/k)+\sum_{j \ioe x/y} \mu(j) (\psi(x/j)-x/j)+xm(x/y)-(\psi(y)-y)M(x/y)-yM(x/y) 
\end{equation*}
mais $xm(x/y)-yM(x/y)=x \left( m(x/y)-\dfrac{M(x/y)}{x/y}  \right)=xm_1(x/y)$.
\end{proof}

\section{$|M(x)|$ pour $x$ entre $10^{16}$ et $2\times 10^{16}$} \label{premier bout}
Dans cet intervalle on tâche d'être le plus précis possible dans nos estimations car c'est ici qu'on trouve les plus grandes valeurs pour nos majorations de $|M(x)|/x$, autrement dit le $1/160\,383$ du théorème \ref{thm A}. 
\subsection{$|N(x)|$ pour $10^{16} \ioe x \ioe 2\times 10^{16}$}
Les estimations $M(x) \ll \sqrt{x} $, $\psi(x)-x \ll \sqrt{x}$ utilisées directement dans l'identité du lemme \ref{identite}  fournissent  $|N(x)| \ll x^{3/4} $ au lieu du $|N(x)| \ll \sqrt{x} \log x $ attendu. Cependant les modèles du lemme \ref{racine} sont numériquement meilleurs que les autres estimations dont on dispose  donc on va utiliser ces modèles le plus possible. 

\begin{prop} \label{premier morceau} On a
\begin{equation*}
    \frac{|N(x)|}{x \log x} \ioe 6.234983 \times 10^{-6} \ioe \frac{1}{160\,385}  \quad (10^{16} \ioe x \ioe 2\times 10^{16})\;. 
\end{equation*}
\end{prop}

On va dans cette partie utiliser le lemme \ref{identite} avec $y'=x/y$ et ainsi calculer la valeur de $m_1(y)$ au lieu de l'estimer :
\begin{equation*}
    \frac{-N(x)}{x}=\sum_{k \ioe x/y} \frac{\Lambda(k)}{k} \frac{M(x/k)}{x/k}\\
    +\sum_{j \ioe y} \frac{\mu(j)}{j} \frac{\psi(x/j)-x/j}{x/j}-\frac{\psi(x/y)-x/y}{x/y}\frac{M(y)}{y}+m_1(y)\;. 
\end{equation*}
Nous profiterons également dans la première des sommes ci-dessus du fait que $\Lambda(1)=0$ pour aller jusqu'à $x=2\times 10^{16}$ alors que le modèle pour $|M|$ s'arrête à $x=10^{16}$.
\begin{lem} \label{avec m1}
Posons $y=10^8$. Si $ 10^{16} \ioe x \ioe 2\times10^{16}$ alors
\begin{equation*}
\frac{1}{\log x} \left( \frac{|\psi(x/y)-x/y|}{x/y}\frac{|M(y)|}{y}+|m_1(y)|\right) \ioe  3.3 \times 10^{-8}\; .
\end{equation*}
\end{lem}
\begin{proof}
On calcule $m_1(y)$ à l'aide du logiciel PARI/GP et on obtient l'encadrement $0<m_1(y)<1.195 \times 10^{-6}$. Puis en utilisant les modèles en racine du lemme \ref{racine} (ce qui est possible car $11 \ioe x/y \ioe 10^{19}$) et la décroissance des fonctions en jeu, on obtient un résultat inférieur à
\begin{equation*} 
\frac{1}{\log 10^{16}} \left( \frac{0.94 \times 0.571}{\sqrt{10^{16}}} +1.195 \times 10^{-6}\right) \ioe 1.46 \times 10^{-10}+ 3.25 \times 10^{-8}\ioe   3.3 \times 10^{-8}\, .
\end{equation*}
\end{proof}

\begin{lem} \label{lambda racine} Pour $11 \ioe X \ioe 10^{19}$ on a
\begin{equation*}
-5.44-0.47 \log X   \ioe   \sum_{k \ioe X} \frac{\Lambda(k)}{\sqrt{k}}-2\sqrt{X} \ioe 0.47 \log X -1.19 \,.
\end{equation*}
\end{lem}

\begin{proof}
La formule sommatoire d'Abel  donne la première égalité qui suit, ensuite on utilise le modèle pour $|\psi(X)-X|$ du lemme \ref{racine} qui n'est valide que pour $X$ plus grand que $11$, pour finir les intégrales sont faciles à calculer,
\begin{align*}
    \sum_{k \ioe X} \frac{\Lambda(k)}{\sqrt{k}}&= \frac{\psi(X)}{\sqrt{X}}+\int_{1}^{X} \psi(t) \frac{t^{-3/2}}{2} \diff{t} \\
    &= \sqrt{X}+ O^{*}(0.94)+ \int_{11}^{X}  \frac{t^{-1/2}}{2} \diff{t}+ O^{*}(0.94)\int_{11}^{X}  \frac{1}{2t} \diff{t}+\int_{1}^{11}\psi(t) \frac{t^{-3/2}}{2} \diff{t} \\
    &=\sqrt{X}+ O^{*}(0.94)+\sqrt{X}-\sqrt{11}+O^*(0.47 \log X)+O^*(1.18)\,.
\end{align*}
\end{proof}

\begin{lem} \label{Lambda debut}
Posons $y=10^8$. Si $ 10^{16} \ioe x \ioe 2\times10^{16}$ alors
\begin{equation*}
\frac{1}{\log x}  \sum_{1<k \ioe x/y} \frac{\Lambda(k)}{k} \frac{|M(x/k)|}{x/k}    
\ioe   3.1002 \times 10^{-6} \, .
\end{equation*}
\end{lem}
\begin{proof}
On utilise le lemme \ref{racine} pour la fonction $|M|$ et le lemme \ref{lambda racine} pour la somme en $\Lambda(k)/\sqrt{k}$ avec $X=x/y$ (ce qui est possible car $11 \ioe x/y \ioe 10^{19}$) puis  la décroissance de $x\mapsto 1/\log x$ ainsi :
\begin{align*}
&\frac{1}{\sqrt{ x }\log x}  \sum_{1<k \ioe x/y} \frac{\Lambda(k)}{\sqrt{k}} \frac{|M(x/k)|}{\sqrt{x/k}}    \\
&   \ioe  \frac{ 1}{ \sqrt{x} \log x} (2  \sqrt{x/y}+0.14\log x -0.14 \log y-1.19)\sup_{y\ioe t \ioe 10^{16} }\frac{|M(t)|}{\sqrt{t}}   \\
&   \ioe   \frac{0.571 \times 2  \sqrt{y}}{ \sqrt{10^{16}} \log 10^{16}}+\frac{0.571 \times 0.14}{\sqrt{10^{16}}}-\frac{0.571 \times (1.19+0.14\log y)}{\sqrt{2 \times 10^{16}} \log(2 \times 10^{16})}                                                      
\\ &\ioe   3.0998 \times 10^{-6}+8 \times 10^{-10}-4 \times 10^{-10} \quad .
\end{align*}
\end{proof}

\begin{lem}\label{avec psi}
Posons $y=10^8$. Si $ 10^{16} \ioe x \ioe 2\times10^{16}$ alors
\begin{equation*}
\frac{1}{\log x}  \sum_{1 \ioe j \ioe y} \frac{\mu^2(j)}{j} \frac{|\psi(x/j)-x/j|}{x/j}\ioe 3.1023 \times 10^{-6}\,.
\end{equation*} 
\end{lem}
\begin{proof}
On calcule $\sum_{1 \ioe j \ioe y} \mu^2(j)/\sqrt{j}$ à l'aide du logiciel PARI/GP et on obtient un résultat inférieur à  $12158.55$. On utilise le lemme \ref{racine} pour la fonction $\psi$ (ce qui est possible car $11 \ioe x/y \ioe x/j \ioe x \ioe 2 \times 10^{16} \ioe 10^{19}$ pour tout $j$ entre $1$ et $y$)
\begin{align*}
&\frac{1}{\sqrt{x}\log x}  \sum_{1 \ioe j \ioe y} \frac{\mu^2(j)}{\sqrt{j}} \frac{|\psi(x/j)-x/j|}{\sqrt{x/j}}\\
\ioe& \frac{0.94}{\sqrt{x} \log x } \sum_{1 \ioe j \ioe y} \frac{\mu^2(j)}{\sqrt{j}} \ioe \frac{0.94}{ \sqrt{10^{16}} \log 10^{16}}  12158.55 \ioe 3.1023 \times 10^{-6}\,.
\end{align*} 
\end{proof}

\subsection{Proximité entre $N(x)$ et $M(x) \log x$, conversion du modèle en racine}
Le modèle du lemme \ref{racine} pour $|M(x)|$ va jusqu'à $10^{16}$ mais en utilisant des méthodes de convolutions intégrales issues de notre thèse \cite{daval2019identites} nous pouvons convertir ce modèle pour $|N(x)|$ en dépassant même un peu $10^{16}$.
\begin{lem} \label{racine N} Pour tous les nombres $x$ vérifiant $1.3 \times 10^9 \ioe x \ioe 6 \times 10^{16}$ on a
 \begin{equation*}
  \Big| \frac{N(x)}{x \log x }- \frac{M(x)}{x} \Big| \ioe \frac{0.227}{\sqrt{x} \log x}    \;.
\end{equation*}
\end{lem}

\begin{cor} \label{conv vers M} On a
\begin{equation*}
    \frac{|M(x)|}{x } \ioe 6.235045 \times 10^{-6} \ioe \frac{1}{160\,383}  \quad (10^{16} \ioe x \ioe 2\times 10^{16}) \;. 
\end{equation*}
\end{cor}

\begin{proof}[Démonstration du lemme \ref{racine N}]

Tout d'abord on a
\begin{equation}
   \frac{N(x)}{x \log x }= \frac{M(x)}{x}-\frac{1}{x \log x} \int_{1}^{x}M(x/t)\frac{1}{t} \diff{t} \quad (x \soe 1) \label{MNint}
\end{equation}
car $\int_{1}^{x}M(t)/t \diff{t}=\sum_{n \ioe x} \mu(n) \log(x/n)=  M(x) \log x -N(x)$.
On considère ensuite la fonction utilisée par Tchebichef dans son étude sur la fonction $\psi$ qui porte depuis son nom
$A(t)=\lfloor t \rfloor- \lfloor t/2 \rfloor- \lfloor t/3 \rfloor- \lfloor t/5 \rfloor+ \lfloor t/30 \rfloor$.
Par l'interversion de somme et d'intégrale suivante (voir notre thèse \cite{daval2019identites} lemme 1 p.5) pour  $\delta : [1, \infty[ \rightarrow \CC$ une fonction localement intégrable on a :
\begin{equation*} 
\int_{1}^{x}M(x/t) \delta(t) \frac{\diff{t}}{t}=\int_{1}^{x} \sum_{n \ioe u} \mu(n)\delta\left(\frac{u}{n}\right) \frac{\diff{u}}{u}\quad (x \soe 1)\;,
\end{equation*}
or par inversion de M\"{o}bius  pour $k\soe 1$ on a $ \sum_{n \ioe u} \mu(n)\lfloor u/(kn) \rfloor=\mathrm{1}_{[k,\infty]}(u)$, donc
\begin{equation*}
 \int_{1}^{x} M(x/t) \frac{A(t)}{t} \diff{t}=  \log x -\log(x/2)-\log(x/3)-\log(x/5)+\log(x/30)=-\log x \quad (x \soe 6) 
\end{equation*}
et ainsi
\begin{equation*}
 \int_{1}^{x} M(x/t) \frac{1}{t} \diff{t}=  \int_{1}^{x} M(x/t) \frac{1-A(t)}{t} \diff{t}-\log x  \quad (x \soe 6) \;.
\end{equation*}
 En utilisant le fait que $A$ est une fonction périodique de période $30$ on vérifie la première des deux propriétés suivantes 
\begin{equation}
0 \ioe   1-A(t) \ioe 1 \quad (t \soe 1) \quad \text { et } \quad 1-A(t)=0   \quad (1 \ioe t < 6) \;.
\end{equation}
Et la positivité permet de calculer facilement l'intégrale suivante
\begin{equation} \label{Mellin}
 \int_{1}^{\infty}|1-A(t)|t^{-1-s} \diff{t}=\frac{1}{s}-\frac{\zeta(s)}{s}(1^{-s}-2^{-s}-3^{-s}-5^{-s}+30^{-s}) \quad (\Re(s)>1, s\neq 1)
\end{equation}
en particulier
\begin{equation*}
 \alpha=\int_{1}^{\infty}|1-A(t)|t^{-3/2} \diff{t}=2-2(1^{-1/2}-2^{-1/2}-3^{-1/2}-5^{-1/2}+30^{-1/2})\zeta(1/2)=0.3962...
\end{equation*}
Ainsi
\begin{align*}
  &\Big| \frac{N(x)}{x \log x }- \frac{M(x)}{x} \Big| \ioe \frac{1}{ \log x}  \int_{6}^{x}\frac{|M(x/t)|}{x/t} \frac{|1-A(t)|}{t^2} \diff{t}+\frac{1}{x}\\
  \ioe& \frac{1}{ \sqrt{x} \log x} \sup_{33\ioe t \ioe x/6 }\frac{|M(t)|}{\sqrt{t}}\int_{6}^{\infty}\frac{|1-A(t)|}{t^{3/2}} \diff{t}+  \frac{1}{\log x}\int_{6}^{x/33}\frac{|M(x/t)|}{x/t} \frac{1}{t^2} \diff{t}+\frac{1}{x}
\end{align*}
or
\begin{equation*}
    \int_{1}^{33}\frac{|M(t)|}{t} \diff{t}=\sum_{n=1}^{32}|M(n)|\log(1+1/n)=4.78...<5
\end{equation*}
donc
\begin{equation}
  \Big| \frac{N(x)}{x \log x }- \frac{M(x)}{x} \Big| \ioe \frac{\alpha \times 0.571 }{\sqrt{x} \log x}  +\frac{1}{x}+ \frac{5}{x \log x}=f(x)
\end{equation}
et $f(x)\sqrt{x} \log x $ est une fonction décroissante qui vérifie $f(1.3 \times 10^9)<0.227$ ce qui donne la majoration annoncée.
\end{proof}
Comparons avec les anciennes méthodes. Schoenfeld dans l'article \cite{SchoenfeldMobius} puis El Marraki dans \cite{MarrakiMobiusIII} passent des majorations de $|N|$ vers des majorations pour $|M|$ à l'aide de l'identité 
\begin{equation*}
\frac{M(x)}{x}-\frac{N(x)}{x \log x}=\frac{M(T)}{x}-\frac{N(T)}{x \log T} + \frac{1}{x}\int_{T}^{x} \frac{N(t)}{ \log^2 t} \frac{\diff{t}}{t}.
\end{equation*}
Avec une majoration pour le logarithme intégral $ \int_{2}^x 1/ \log t \diff{t} \ioe x / \log x +O(1/\log^2x) $ on obtient
\begin{equation*}\Big|\frac{M(x)}{x}   -\frac{N(x)}{x \log x }\Big| \ioe \frac{1}{ \log x}\sup_{T < u < x} \frac{|N(u)|}{u \log u} +O( \log^{-2} x) 
 \end{equation*}
qui aurait fournit au mieux
\begin{equation*}
    \frac{|M(x)|}{x} \ioe (1+1/\log(10^{16}))\frac{1}{160\,385} \ioe \frac{1}{156\,146}\,.
\end{equation*}

\section{$|N(x)|$ et $|M(x)|$ pour $x$ entre $2\times 10^{16}$ et $ 10^{19}$}
Nous allons montrer que $|M(u)|/u \ioe 1/160\,383$ jusqu'à $10^{19}$ en découpant l'analyse.
\begin{lem} \label{N M 1}
Si $|M(u)|/u \ioe 1/160\,383$ pour tout $u$ vérifiant $(0.571 \times 160\,383)^2\ioe u \ioe x/6$ alors
 \begin{equation*}
  \Big| \frac{N(x)}{x \log x }- \frac{M(x)}{x} \Big| \ioe  \frac{0.08}{160\,383} \frac{1}{\log x}+\frac{104\,586}{x \log x} + \frac{1}{x}  \;.
\end{equation*}
Si de plus $x\soe 2 \times 10^{16}$ alors
\begin{equation*}
  \Big| \frac{N(x)}{x \log x }- \frac{M(x)}{x} \Big| \ioe 1.33 \times 10^{-8}  \;.
\end{equation*}
\end{lem}

Nous utiliserons désormais la formule du lemme \ref{identite} qui comporte le terme $m_1(x/y)$ donc nous avons besoin d'estimer cette fonction.  Remarquons que pour $y=10^8$ dans la section \ref{premier bout} on a obtenu par calcul direct $m_1(y) \approx 0.0119/\sqrt{y}$.
\begin{lem} [Modèles en racine carrée]\label{racine m1} On a 
\begin{align*}
|m_1(X)| &\ioe 0.129/\sqrt{X} \quad    (5\times 10^6 \ioe X\ioe 10^{16})\,, \\
|m_1(X)| &\ioe  0.129/\sqrt{A}   \quad (X \soe A=2.18 \times 10^{12})\,.
\end{align*} 
\begin{proof}
Voir la prépublication \cite[théorème 2 et corollaire 1]{daval2020conversions}. 
\end{proof}
\end{lem}
\begin{lem} \label{racine mu carre} On a 
\begin{equation*}
   \frac{1}{\sqrt{x}} \sum_{k \ioe x/y} \frac{\mu^2(k)}{\sqrt{k}} \ioe  \frac{12}{ \pi^2 \sqrt{y}}+\frac{0.17 \log (x/y)}{\sqrt{x}}  \quad  (x/y \soe 10^4) \,.
\end{equation*}
\end{lem}
\begin{lem} \label{estimation simple} On a
\begin{equation*}  \sum_{k \ioe y} \frac{\Lambda(k)}{k} \ioe \log y  \quad (y \soe 1) \quad \text{ et } \quad  \sum_{k \ioe y} \frac{\Lambda(k)}{k} \ioe \log y -0.5 \quad (y \soe 300) \,.
\end{equation*}
\end{lem}
\begin{proof}
Notons $R(y)=\sum_{k \ioe y} \Lambda(k)/k-\log x+\gamma$ d'après l'article \cite{ramare2013explicit} pour $1 \ioe y \ioe 10^{10}$ on a $|R(y)|\ioe 1.31/ \sqrt{y}$ et pour $y \soe 10^{10}$ on a $|R(y)| \ioe 0.0067/\log y$. Voir également le corrigendum \cite{ramarecorrigendum}.
\end{proof}

\begin{lem} \label{dyadique}
Posons  $X=a\times 10^{16}$ avec $0<a\ioe 500$ et $y=\sqrt{X}$. Si $X \ioe x \ioe 2X$ alors
 \begin{align*}
     \frac{|N(x)|}{x \log x} \ioe \frac{1}{\log X} \left( \frac{0.94 \times 0.571}{\sqrt{X}}+\frac{0.129}{\sqrt{ X/y}}\right)+ \frac{0.94}{  \log X}\frac{12}{ \pi^2 \sqrt{y}}+\frac{0.17 \times 0.94 }{\sqrt{X}} \\
     +\frac{\log(2a-1)}{  \log X}   \sup_{10^{16}\ioe t\ioe 2X}\frac{|M(t)|}{t}+  \frac{  0.571}{ \sqrt{X} \log X} (2  \sqrt{y}+0.14\log y-1.19 ) \,.
 \end{align*}
\end{lem}
On applique successivement le lemme précédent avec $a=2, 4, 8, \dots,  256$ et pour finir $500$,  en majorant à chaque fois $\sup_{10^{16}\ioe t\ioe 2X}|M(t)|/t$ par $1/160\,383$, puis on utilise le lemme \ref{N M 1} pour passer de $|N|$ à $|M|$  et on vérifie à chaque itération que  $\sup_{ X\ioe x\ioe  2X}   |M(x)|/x \ioe 1/160\,383$. On résume ces étapes avec le tableau suivant. 
\[
\begin{array}{|c|c|c|c|c|c|c|c|c|c|} \hline
a  \text{ et } X=10^{16}a   & 2     & 4        & 8                & 16      & 32     & 64      &128 &256  &500   \\ \hline
10^6 \sup_{ X\ioe x\ioe  2X} \frac{|N(x)|}{x \log x} \ioe \phantom{\dfrac{1}{2}}  & 5.60   & 4.79     &    4.13         &  3.59    & 3.16 & 2.82 &  2.56   & 2.35 & 2.19 \\  \hline 
10^6 \sup_{ X\ioe x\ioe  2X}   \frac{|M(x)|}{x } \ioe  \phantom{\dfrac{1}{2}}        &  5.61      & 4.80      &  4.14   & 3.60  & 3.18 &    2.84   & 2.57 & 2.36& 2.20 \\  \hline
\sup_{ X\ioe x\ioe  2X}   \frac{|M(x)|}{x } \ioe \frac{1}{160\,383}     \phantom{\dfrac{1}{2}}     &  \text{oui}     & \text{oui}      &  \text{oui}  & \text{oui} &\text{oui}  & \text{oui} &  \text{oui}   & \text{oui} & \text{oui} \\  \hline
\end{array}
\]
\begin{proof}[Démonstration du lemme \ref{dyadique}]
On a $y=\sqrt{X}=X/y$ donc $\sqrt{2} \times 10^{8} \ioe y \ioe 10^{9.5}$ donc on peut utiliser les modèles en racine (lemme \ref{racine} pour $\psi$ et $M$, lemme \ref{racine m1} pour $m_1$) ainsi on majore deux des quatre termes présents dans le lemme \ref{identite} :
\begin{equation*}
\frac{1}{\log x} \left( \frac{|\psi(y)-y|}{y}\frac{|M(x/y)|}{x/y}+|m_1(x/y)|\right) \\
\ioe  \frac{1}{\log X} \left( \frac{0.94 \times 0.571}{\sqrt{X}}+\frac{0.129}{\sqrt{ X/y}}\right)\,.
\end{equation*}
De même on peut utiliser le lemme \ref{racine mu carre} pour la somme en $\mu^2(j)/\sqrt{j}$ et le modèle en racine pour $\psi$ précédemment invoqué :
\begin{align*}
\frac{1}{\log x}  \sum_{1 \ioe j \ioe x/y} \frac{\mu^2(j)}{j} \frac{|\psi(x/j)-x/j|}{x/j}
\ioe \frac{0.94}{\sqrt{x} \log x } \sum_{1 \ioe j \ioe x/y} \frac{\mu^2(j)}{\sqrt{j}}
\ioe \frac{0.94}{  \log X}\frac{12}{ \pi^2 \sqrt{y}}+\frac{0.17 \times 0.94 }{\sqrt{X}} \,.
\end{align*} 
Le dernier terme est le seul où le fait que $x \ioe 2X$ soit utilisé. Séparons les $x/k$ éventuellement plus grands que $10^{16}$ à l'aide de $a=X/10^{16}$, par le lemme \ref{estimation simple} pour la somme en  $\Lambda(k)/k$ on arrive à :
\begin{align*}
&\frac{1}{\log x}  \sum_{k \ioe y} \frac{\Lambda(k)}{k} \frac{|M(x/k)|}{x/k}=\frac{1}{\log x}  \sum_{1 \ioe k \ioe 2a-1} \frac{\Lambda(k)}{k} \frac{|M(x/k)|}{x/k}+\frac{1}{\sqrt{ x }\log x}  \sum_{2a \ioe k \ioe y} \frac{\Lambda(k)}{\sqrt{k}} \frac{|M(x/k)|}{\sqrt{x/k}}    \\
&   \ioe \frac{1}{  \log X}  \sum_{k \ioe 2a-1} \frac{\Lambda(k)}{k} \sup_{10^{16}\ioe t\ioe 2X}\frac{|M(t)|}{t}+  \frac{  0.571}{ \sqrt{X} \log X} (2  \sqrt{y}+0.14\log y-1.19 ) \,.
\end{align*}
\end{proof}
\section{$|N(x)|$ et $|M(x)|$ pour  $x \soe  10^{19}$}

Dans cette section nous pouvons utiliser des majorations plus grossières que pour $x$ entre $10^{16}$ et $10^{19}$. 

Nous présentons les deux lemmes suivants pour $y=y_1=(0.571\times4345)^2$ mais comme le montrent leurs démonstrations, ils sont vrais pour toutes valeurs de $y$ telles que $y^2 \ioe 10^{19}$ et telles que $|M(x)|/x \ioe 0.571 /\sqrt{y}$ pour tout $x \soe y$.

\begin{lem} \label{N apres 10 puissance 19}
Posons $y=y_1=(0.571\times4345)^2$, $Y=\exp(40)$ et $A=2.18 \times 10^{12}$. Si $x \soe  10^{19} \soe y^2$  on a
\begin{align*}
\frac{|N(x)|}{x \log x}\ioe  &\frac{1}{\log (  10^{19})} \left( \frac{0.94}{\sqrt{y}}\times \frac{0.571}{\sqrt{ y}}+\frac{0.129}{\sqrt{\min(10^{19}/y, A)}}\right)  
+ \frac{0.571}{\sqrt{y}}\frac{\log y-0.5}{ \log (10^{19})} \\
+&\frac{6}{\pi^2}  \sup_{t \soe Y}\frac{|\psi(t)-t|}{t}+\frac{0.94}{\log 10^{19}}  \left( \frac{12}{\pi^2\sqrt{y}}+\frac{0.35 \log(Y/y)}{\sqrt{10^{19}}} \right)  \,.
\end{align*}
\end{lem}

\begin{lem} \label{N vers M plus loin}
Posons $y=y_1=(0.571\times4345)^2$ on a $33\ioe y \ioe 10^{16}$ et pour tout $x \soe y$ on a
\begin{equation*}
  \frac{|M(x)|}{x  } \ioe  \frac{|N(x)|}{x \log x }+   \frac{0.08}{\log x} \frac{0.571}{\sqrt{y}}+ \frac{5+0.571\times2\sqrt{y}}{x \log x} + \frac{ 1}{x}  \;.
\end{equation*}
\end{lem}
Avec $y=y_1=(0.571\times4345)^2$ qui correspond à la majoration $|M(x)|/x \ioe 1/4345$ pour tout $x\soe y_1$, on obtient que pour tout $x\soe 10^{19}$, $|N(x)|/(x \log x) \ioe 1/ 11086 $ et donc $|M(x)|/x \ioe 1/11035$. Ce résultat et le modèle $|M(x)| \ioe 0.571 \sqrt{x}$ du lemme \ref{racine} fournissent un nouvel $y=y_2=(0.571 \times 11035)^2$ admissible et on recommence avec cette nouvelle majoration pour $|M|$. Pour pouvoir itérer il faut que $y^2 \ioe 10^{19}$ donc on pose $L^*$ qui est tel que $\sqrt{10^{19}}=(0.571 \times L^*)^2$. On a $L^* \approx 98\,483$ et on résume nos itérations par le tableau ci-dessous.
\[
\begin{array}{|c|c|c|c|c|c|} \hline
y=(0.571 \times L)^2 \text{ avec } L=  &   4345   &   11035     &    25266            &    53119  &  L^*   \\ \hline
\sup_{x \soe 10^{19}} \frac{|N(x)|}{x \log x} \ioe \phantom{\dfrac{1}{2}}  & 1/11086   & 1/25372     &    1/53324         &  1/104069   & 1/180799  \\  \hline 
  \sup_{x \soe 10^{19}} \frac{|M(x)|}{x } \ioe \phantom{\dfrac{1}{2}}       &  1/11035      & 1/25266     &  1/53119  & 1/103697  & 1/180194  \\  \hline
\end{array}
\]
Nous allons utiliser les estimations explicites suivantes.
\begin{lem} \label{estimations3} On a
\begin{equation*}
    |\psi(X)-X|\ioe 8.6386 \times 10^{-8}X  \quad (X \soe \exp(40)\approx 2.35 \times 10^{17}) \;,
\end{equation*}
\begin{equation*}
    \sum_{n \ioe X} \frac{\mu^2(n)}{n} \ioe \frac{6}{\pi^2} \log (7X)  \quad (X \soe 1)  \label{mu2} \;,
\end{equation*}
 \begin{equation*}
     \sum_{x/Y<j \ioe x/y} \frac{\mu^2(j)}{j} \frac{1}{\sqrt{x/j}} \ioe \frac{6}{\pi^2}\left(\frac{2}{\sqrt{y}}-\frac{2}{\sqrt{Y}}\right)+ \frac{1.31 \sqrt{Y}}{x}+\frac{0.35 \log(Y/y)}{\sqrt{x}} \quad (x \soe Y \soe y \soe 1) \;.
\end{equation*}
\end{lem}

\begin{proof}[Démonstration du lemme \ref{N apres 10 puissance 19}]
Nous allons majorer l'identité du lemme \ref{identite}. Nous pouvons utiliser les modèles en racine carrée des lemmes \ref{racine} et \ref{racine m1} pour $\psi$ et $m_1$. On réécrit la majoration \eqref{CDM} $|M(u)|/u \ioe 1/4345 $ pour tout $u\soe y=(0.571 \times 4345)^2 $ en   $|M(u)|/u \ioe 0.571/\sqrt{y} $ pour tout $u\soe y$.
Puisque $y \ioe x^2$ on  a $x/y \soe y$ donc on a
\begin{equation*}
\frac{1}{\log x} \left( \frac{|\psi(y)-y|}{y}\frac{|M(x/y)|}{x/y}+|m_1(x/y)|\right) 
\ioe  \frac{1}{\log (  10^{19})} \left( \frac{0.94}{\sqrt{y}}\times \frac{0.571}{\sqrt{ y}}+\frac{0.129}{\sqrt{\min(10^{19}/y, A)}}\right) . 
\end{equation*}

Passons à la contribution la plus importante. Remarquons qu'en suivant la méthode de Schoenfeld de l'article \cite{SchoenfeldMobius} cette valeur serait doublée.

\begin{equation*}
\frac{1}{\log x}  \sum_{k \ioe y} \frac{\Lambda(k)}{k} \frac{|M(x/k)|}{x/k} \ioe  \sup_{t \soe x/y}\frac{|M(t)|}{t} \frac{\log y-0.5}{ \log x} \ioe \frac{0.571}{\sqrt{y}}\frac{\log y-0.5}{ \log (10^{19})} \, .
\end{equation*} 

On va découper la dernière somme à majorer en deux et utiliser des majorations pour $\psi$ différentes. 
\begin{equation*}
\frac{1}{\log x}  \sum_{j \ioe x/Y} \frac{\mu^2(j)}{j} \frac{|\psi(x/j)-x/j|}{x/j} \ioe  \frac{6}{\pi^2}  \sup_{t \soe Y}\frac{|\psi(t)-t|}{t}  \,.
\end{equation*} On utilise  le lemme \ref{estimations3}, et le fait que $-\frac{12}{\pi^2\sqrt{Y}}+ \frac{1.31 \sqrt{Y}}{x }<0$ car $Y<0.92X$, pour obtenir 
\begin{align*}
\frac{1}{\log x}  \sum_{x/Y<j \ioe x/y} \frac{\mu^2(j)}{j} \frac{|\psi(x/j)-x/j|}{x/j} 
\ioe& \frac{1}{\log x} \sup_{ y\ioe t \ioe Y}\frac{|\psi(t)-t|}{\sqrt{t}} \left( \frac{12}{\pi^2\sqrt{y}}+\frac{0.35 \log(Y/y)}{\sqrt{x}} \right)  \\
\ioe &\frac{0.94}{\log 10^{19}}  \left( \frac{12}{\pi^2\sqrt{y}}+\frac{0.35 \log(Y/y)}{\sqrt{10^{19}}} \right) \,. 
\end{align*} 
\end{proof}

\begin{proof}[Démonstration du lemme \ref{estimations3}]
L'estimation pour $\psi$ provient de l'article \cite{faber2015new} rectifié par le corrigendum \cite{faber2018corrigendum}.

L'article \cite[lemma 3.4]{ramare2016explicit} donne $0.578 \ioe \sum_{n \ioe X} \mu^2(n)/n-6/\pi^2 \log X \ioe 1.166$ pour tout nombre réel $X \soe 1$, et on a $\exp(1.166/(6/\pi^2))<7$.

On utilise l'article  \cite[lemma 4.7 p.~370]{ramare2013explicit} avec la fonction $f(x)=\sqrt{x}$ qui est bien positive, croissante et continûment dérivable :
\begin{equation*}
     \sum_{x/Y<j \ioe x/y} \mu^2(j) \sqrt{x/j} \ioe 1.31 \sqrt{Y}+\frac{6x}{\pi^2} \int_{y}^{Y} t^{-1.5} \diff{t}+0.35\sqrt{x} \int_{y}^{Y} t^{-1} \diff{t} \,.
\end{equation*}
Pour finir  les intégrales ci-dessus se calculent de manière exacte. 
\end{proof}

\section{Reste des démonstrations des sections précédentes}

\begin{proof}[Démonstration du lemme \ref{racine}]
L'estimation sur $M$ provient de l'article \cite{hurst2018computations} et l'estimation pour $\psi$ provient de l'article \cite[theorem 2]{buthe2018analytic}.
\end{proof}

\begin{proof}[Démonstration du lemme \ref{racine mu carre}]
Notons $Q(x)=\sum_{k \ioe x} \mu^2(k)$ on a $Q(x) \ioe (6/\pi)^2x+0.1333 \sqrt{x}$ pour $x\soe 1664$ et $Q(x) \ioe (6/\pi)^2x+0.5 \sqrt{x}$  pour $x\soe 1$. On note $X=x/y$
 La formule sommatoire d'Abel  donne la première égalité qui suit
\begin{align*}
    \sum_{k \ioe X} \frac{\mu^2(k)}{\sqrt{k}}=& \frac{Q(X)}{\sqrt{X}}+\int_{1}^{X} Q(t) \frac{t^{-3/2}}{2} \diff{t}\\
    \ioe  &\frac{6}{\pi^2}\sqrt{X}+ \frac{6}{\pi^2}\int_{1}^{X}  \frac{t^{-1/2}}{2} \diff{t}+0.1333(1+0.5\log X)+R 
\end{align*}
\begin{equation*}
    \text{où } R=\int_{1}^{1664}|Q(t)-(6/\pi)^2t-0.1333 \sqrt{t}|\frac{t^{-3/2}}{2} \ioe 0.5(0.5-0.1333)\log 1664  \ioe 1.36 \,.
\end{equation*}
Et pour finir $-(6/\pi^2)\sqrt{1}+0.1333+R \ioe 0.1 \log (x/y)$ si $x/y \soe 10^4$.
\end{proof}

\begin{proof}[Démonstration des lemmes \ref{N M 1} et \ref{N vers M plus loin}]
On suit le même principe que dans la démonstration du lemme \ref{racine N}. On passe à la limite quand $s \to 1$ dans \eqref{Mellin} en utilisant $\zeta(s)=1/(s-1)+O(1)$ (il est possible d'intervertir car l'intégrale converge absolument pour $\Re(s)>0$).
Posons 
\begin{equation*}
 \beta=\int_{1}^{\infty}|1-A(t)|t^{-2} \diff{t}=1+\frac{\log 1}{1}-\frac{\log 2}{2}-\frac{\log 3}{3}-\frac{\log 5}{5}+\frac{\log 30}{30}=0.07870...
\end{equation*}
alors
\begin{equation}
 | \frac{1}{x} \int_{1}^{x} M(x/t) \frac{1}{t} \diff{t}| \ioe   \beta \sup_{T \ioe u \ioe x/6} \frac{|M(u)|}{u} + \frac{1}{x} \int_{1}^{T} |M(t)| \frac{1}{t} \diff{t}+ \frac{\log x}{x}
\end{equation}
et si $33\ioe T\ioe 10^{16}$ on a
\begin{equation}
 | \frac{1}{x} \int_{1}^{x} M(x/t) \frac{1}{t} \diff{t}| \ioe   \beta \sup_{T \ioe u \ioe x/6} \frac{|M(u)|}{u} + \frac{5+0.571\times2\sqrt{T}}{x} + \frac{\log x}{x} \, .
\end{equation}
Pour conclure on choisit $T=(0.571 \times 160\,383)^2$ qui est le plus petit nombre tel que $|M(x)|/x \ioe 0.571/\sqrt{x} \ioe 1/160\,383$ pour tout $x$ entre $T$ et $10^{16}$ suivant le modèle en racine du lemme \ref{racine}.
\end{proof}

\begin{proof}[Démonstration de la proposition \ref{premier morceau}]
On majore l'identité du lemme \ref{identite}, pour cela on reprend les lemmes \ref{avec m1}, \ref{Lambda debut} et \ref{avec psi} sans utiliser la valeur arrondie de chaque énoncé qui n'est présente qu'à titre indicatif.
\end{proof}

\begin{proof}[Démonstration du théorème \ref{thm A}]
Pour démontrer le résultat pour $x\soe 10^{16}$ on recolle le corollaire \ref{conv vers M} et le tableau situé après après le lemme \ref{N vers M plus loin}. Pour finir on utilise le modèle du lemme \ref{racine} $|M(x)| \ioe 0.571 \sqrt{x}$ qui est valable entre $33$ et $10^{16}$ :
\begin{equation*}
\frac{|M(x)|}{x} \ioe  \frac{1} {160\,383}  \quad (x \soe  10^{16}) 
 \text{ , }    \frac{1} {160\,383} \ioe \frac{0.571}{\sqrt{x}}
\quad (T \ioe x \ioe 10^{16}) 
\end{equation*}
où $T=(0.571 \times 160\,383)^2 \soe 8.3867 \times 10^{9}$.
\end{proof}

\begin{rem}
 On pourrait compléter le résultat démontré ci-dessus en cherchant informatiquement le plus petit rang avant $8.4 \times 10^9$ pour lequel $|M(x)|/x \ioe 160\,383$.
\end{rem}

\bibliographystyle{alpha}

\bibliography{bibarticle}

\end{document}